\documentclass{elsarticle}

\biboptions{sort&compress,comma,round}
\usepackage{lineno,url}
\usepackage{amsmath,amsfonts,amssymb}
\usepackage{amssymb,amsmath,amsthm}
\usepackage{overpic}
\usepackage{contour}\contourlength{0.7pt}
\usepackage{tikz}\usetikzlibrary{shapes,arrows}
\usepackage{versions}
\usepackage{natbib}
\usepackage{algorithm, algorithmic}
\usepackage{subcaption}

\newtheoremstyle{mytheorem}{5pt plus 5pt minus 3pt}{4pt plus 3pt minus 1.5pt}
	{\itshape}{}{\bfseries}{.}{1ex plus 1ex minus .5ex}{}
\newtheoremstyle{mydef}{5pt plus 5pt minus 3pt}{4pt plus 3pt minus 1.5pt}
	{}{0pt}{\bfseries}{.}{1ex plus 1ex minus .5ex}{}
\newtheoremstyle{myremark}{5pt plus 5pt minus 3pt}{4pt plus 3pt minus 1.5pt}
	{}{0pt}{\itshape}{.}{1ex plus 1ex minus .5ex}{}
\theoremstyle{mytheorem}

\newtheorem{rem}{Remark}
\theoremstyle{mydef}

\theoremstyle{myremark}

\newcommand{\bi}{\mathbf i}
\newcommand{\bj}{\mathbf j}

\newcommand{\bb}{\mathbf b}

\newcommand{\bg}{\mathbf g}

\newcommand{\de}{\mathrm{d}}

\newcommand{\bx}{\mathbf x}
\newcommand{\bu}{\mathbf u}

\newcommand{\Om}{\Omega}
\newcommand{\Omt}{\widetilde{\Omega}}

\def\Transp{^{\text{\sf\sc t}}}
\newcommand{\RR}{\mathbb R}
\newcommand{\ZZ}{\mathbb Z}
\newcommand{\bxh}{\hat{\bx}}
\newcommand{\xh}{\hat{x}}
\newcommand{\uh}{\hat{u}}
\newcommand{\Omh}{\hat{\Omega}}
\newcommand{\Bh}{\hat{B}}
\newcommand{\Sh}{\hat{S}}
\newcommand{\noe}[2]{N^\textnormal{#1}_{#2}}
\newcommand{\nh}{\hat{\nabla}}
\newcommand{\bK}{\mathbf K}
\newcommand{\bM}{\mathbf M}

\def\Acaption#1#2{\caption{#2}\vspace*{-#1}}

\definecolor{Mgreen}{RGB}{34,139,34}
\definecolor{blau}{rgb}{0.15,0.2,0.5}
\definecolor{gray}{rgb}{0.5,0.5,0.5}
\definecolor{drot}{rgb}{0.7,0,0.1}
\definecolor{gelb}{rgb}{.55,.40,.1}
\definecolor{magenta}{rgb}{1.,0.,1.}
\definecolor{cyan}{rgb}{0.,1.,1.}
\definecolor{green}{rgb}{0.,1.,0.}
\definecolor{Morange}{rgb}{1.,0.5,0.}

\journal{CMAME}

\begin{document}

\begin{frontmatter}
\title{Efficient mass and stiffness matrix assembly via weighted Gaussian quadrature rules for B-splines}

\author[BCAM]{Michael Barto\v{n}\corref{cor1}}
\ead{mbarton@bcamath.org}
\author[Curtin,cic]{Vladimir Puzyrev}
\ead{Vladimir.Puzyrev@Curtin.edu.au}
\author[Curtin,cic]{Quanling Deng}
\ead{Quanling.Deng@Curtin.edu.au}
\author[Curtin,CSIRO,cic]{Victor Calo}
\ead{Victor.Calo@Curtin.edu.au}

\cortext[cor1]{Corresponding author}

\address[BCAM]{BCAM -- Basque Center for Applied Mathematics, Alameda de Mazarredo 14,\\ 48009 Bilbao, Basque Country, Spain}
\address[Curtin]{Department of Applied Geology, Western Australian School of Mines, Faculty of Science and Engineering, Curtin University,
Kent Street, Bentley, Perth, WA 6102, Australia}
\address[CSIRO]{Mineral Resources, Commonwealth Scientific and Industrial Research Organization (CSIRO), Kensington, Perth, WA 6152, Australia}
\address[cic]{Curtin Institute for Computation,
Curtin University, Kent Street, \\Bentley, Perth, WA 6102, Australia}

\begin{abstract}
Calabr{\`o} et al.~\cite{Sangali-2016-Weighted} changed the paradigm of the mass and stiffness computation from the traditional element-wise assembly to a row-wise concept, showing that the latter one offers integration that may be orders of magnitude faster. Considering a B-spline basis function as a non-negative measure, each mass matrix row is integrated by its own quadrature rule with respect to that measure. Each rule is easy to compute as it leads to a linear system of equations, however, the quadrature rules are of the Newton-Cotes type, that is, they require a number of quadrature points that is equal to the dimension of the spline space. In this work, we propose weighted quadrature rules of Gaussian type which require the minimum number of quadrature points while guaranteeing exactness of integration with respect to the weight function. The weighted Gaussian rules arise as solutions of non-linear systems of equations. We derive rules for the mass and stiffness matrices for uniform $C^1$ quadratic and $C^2$ cubic isogeometric discretizations.
Our rules further reduce the number of quadrature points  by a factor of $(\frac{p+1}{2p+1})^d$ when compared to \cite{Sangali-2016-Weighted}, $p$ being the polynomial degree and $d$ the dimension of the problem, 
and consequently reduce the computational cost of the mass and stiffness matrix assembly by a similar factor. 
\end{abstract}


\end{frontmatter}

\section{Introduction}\label{intro}
Systems of Partial Differential Equations (PDEs) describe 
many relevant physical processes.
These physical phenomena are traditionally modeled using finite element analysis (FEA) and Isogeometric analysis (IGA) \cite{Cottrell-2009}. 
%
Numerical integration is a fundamental step in the assembly  process of the algebraic systems that result from the FEA and IGA discretizations. Depending on the type of the governing PDE, the matrix assembly requires integration of products of basis functions (mass matrix) and/or their derivatives (stiffness matrix) which, in the context of IGA, requires integration of spline spaces of a certain structure \cite{StiffnessMatrix-2016}. Efficient quadrature rules play a key role in this process, as they are cheap and elegant tools to exactly integrate the spline space under consideration \citep{Hughes-2010,Calabro-2013-QuadratureNURBS,Oliveira-2009-Weighted,Gautschi-2000-GaussWeight}.

Traditionally in the tensor product based IGA, mass and stiffness matrix assembly is performed element-wise, using a corresponding univariate quadrature rule on each element in each parameter direction. While using standard polynomial Gauss quadrature on each element is a ussual approach in many FEA and IGA codes, Gaussian quadrature for splines offer a lot cheaper alternative as the higher continuity between elements signifies that fewer Gaussian quadrature points are needed \cite{Schoenberg-1958,Micchelli-1977}. 
For example, for the $C^1$ quadratic spline space, the mass matrix contains terms that belong to a quartic $C^1$ space. Gaussian quadrature for this space requires asymptotically, i.e. for a large number of elements, only one and half quadrature points per element in contrast to the polynomial Gauss rule that requires three quadrature points per element \cite{StiffnessMatrix-2016}.

In the IGA community, several recent papers focused on the development of  efficient quadrature rules \cite{Calabro-2013-QuadratureNURBS, Hughes-2010,Hughes-2012,Adam-2015-Selective}. Hughes et al. \cite{Hughes-2010} introduced efficient rules 
that are exact over the whole real line (infinite domain). 
For finite domains, one may introduce additional quadrature points \cite{Hughes-2012} which make the rule non-Gaussian (slightly sub-optimal in terms of quadrature points), but more importantly, it yields quadrature weights that can be negative, unlike Gaussian quadratures.

To overcome these drawbacks, alternative quadrature schemes were introduced in recent years \cite{Bremer-2010-Nonlinear,Schillinger-2014-Reduced,Hiemstra-2017-Optimal,Johannessen-2017-NumOpt,StiffnessMatrix-2016,StiffnessMatrixEven-2016}. In these studies, the sought Gaussian quadrature rule is represented as a root of a piece-wise polynomial system, where the system expresses the exactness of the rule when applied to a basis of the underlying spline space.
In general, the system is highly non-linear and computing a root numerically using, e.g., Newton-Raphson may not always converge \cite{Johannessen-2017-NumOpt}, unless a very good initial guess is known. 
Finding a good initial guess, however, is possible only for specific target spaces using the local structure of the B-spline basis \cite{Adam-2015-Selective,Hiemstra-2017-Optimal}.

An alternative good initial guess has been proposed by using the continuity argument between a spline space and its Gaussian quadrature \cite{Homotopy-2016}. A Gaussian quadrature rule for a desired spline space is derived from a known Gaussian quadrature (e.g., a union of polynomial Gauss rules) by continuously modifying the source knot vector into the target one.
The process is effective for arbitrary knot vectors (including non-uniform spacing and arbitrary continuity) and polynomial degrees.
The sought quadrature rule corresponds to a zero of a piece-wise polynomial system and is traced numerically as the knot vector changes from the source to the target one. This homotopic continuation approach is used for spline spaces of various degrees and continuities \cite{StiffnessMatrix-2016,StiffnessMatrixEven-2016}, showing also the numerical evidence that, for uniform knot vectors, the rules over finite domains quickly converge to the half-point rules of Hughes et al. over infinite domains \cite{Hughes-2010}.

A quadrature-free approach to assemble mass and stiffness matrices uses the observation that exact integration is not required to achieve the optimal convergence rate of the solution. Therefore, the integrals arising from the geometry factor can be approximated by the integrals of the B-spline basis functions, which are precomputed and stored in a look-up table \cite{Mantzaflaris-2012-Exploring,Mantzaflaris-2015-Integration}. 
Another efficient alternative is the variational collocation method \cite{Gomez-2016-Variational}. In \cite{Gomez-2016-Variational}, the authors 
prove the existence of Cauchy-Galerkin collocation points, 
that is, points in which the collocated solution reproduces the Galerkin counterpart. Therefore, such a method has a great advantage as it possess the exactness of Galerkin solution for the cost of collocation. However, a stable computational framework to efficiently determine the Cauchy-Galerkin points is an open challenge, particularly for spline spaces of various continuities and non-uniform knot vectors \citep{Gomez-2017-Variational2}.

Recently, Calabr{\`o} et al. \cite{Sangali-2016-Weighted} have changed the paradigm of the mass and stiffness computation from the traditional element-wise assembly to a row-wise concept. When building the mass matrix, one B-spline basis function of the scalar product is considered as a positive measure (i.e., a weight function), and a weighted quadrature with respect to that weight is computed for each matrix row. Such an approach brings significant 
computational savings compared to the traditional approaches that use Gaussian or semi-Gaussian  element-wise assembly \cite{Schillinger-2014-Reduced,Hiemstra-2017-Optimal,Johannessen-2017-NumOpt,StiffnessMatrix-2016,StiffnessMatrixEven-2016}, because in these integration schemes, the number of quadrature points contains a $p^d$ term, $p$ being the polynomial degree and $d$ the dimension. In contrast, \cite{Sangali-2016-Weighted} requires in each parameter direction in the limit only two points per element, regardless the polynomial degree. For each weight (mass matrix row), its specific weighted quadrature is computed by solving a linear system. These rules, however, are quadratures of the Newton-Cotes type, that is, they require the same number of quadrature points as the spline basis functions involved.

In this work, we propose weighted quadrature rules of Gaussian type which require the minimum number of quadrature points while guaranteeing exactness of integration with respect to the weight function. Our rules further reduce the number of quadrature points of \cite{Sangali-2016-Weighted} from two to one per element. We derive rules for $C^1$ quadratic and $C^2$ cubic isogeometric discretizations over uniform knot vectors. The rules arise as solutions of the resulting non-linear systems of equations. The solution to the quadratic case can be expressed symbolically. For cubic splines, the rules are derived using the Newton-Raphson method. 

The rest of the paper is organized as follows. Section~\ref{sec:Prelim} formulates the model problem and its isogeometric discretization. Then we derive the weighted Gaussian quadrature rules for quadratic and cubic spline spaces in Section~\ref{sec:WQ} and show numerical experiments that validate our theoretical results in Section~\ref{sec:NUMER}. We conclude the paper and indicate directions for future research in Section~\ref{sec:con}. 



\section{Preliminaries}\label{sec:Prelim}

We set up our model problem, derive its isogeometric discretization, and define spline spaces that we need for the approximation of the solution. 

\subsection{Model problem}\label{ssec:MP}
We consider the Poisson problem
\begin{equation}\label{eq:Poisson}
\begin{cases}
\begin{aligned}
-\Delta u & = f	   \qquad\text{in}\> \Omega,\\
           u & = 0		\qquad\text{on}\> \partial\Omega
\end{aligned}
\end{cases}
\end{equation}
as a model problem with a physical domain $\Omega \subset \RR^d$, $d=1,2,3$, with open bounded Lipschitz boundary $\partial\Omega$.  Let $V = H_0^1(\Omega)$, the variational  form of \eqref{eq:Poisson} is to find $u\in V$ such that
\begin{equation}\label{eq:Variational}
a(u,v) = l(v) \quad \forall v \in V,
\end{equation}
where the bilinear and linear forms are 
\begin{equation}\label{eq:Forms}
a(u,v) = \int_{\Omega} \nabla u\Transp(\bx) \cdot \nabla v (\bx) \,\de \bx \quad \textnormal{and} \quad
l(v) = \int_{\Omega} f(\bx) v(\bx) \,\de \bx,  
\end{equation}
respectively. 

\subsection{Isogeometric discretization}\label{ssec:IGA}

In the isogeometric framework, $\Omega$ is described with the same basis as the solution space, that is, it is an isoparametric discretization. In this case, we commonly use B-spline or NURBS parameterizations.
Consider a single patch geometry map $G: \Omh \rightarrow \Om$, where the parameter domain $\Omh$ is a unit box in $\RR^d$, $\Omh = [0,1]^d$. The mapping $G$ maps any point $\bxh \in \Omh$  from the parameter domain to the physical domain via 
\begin{equation}\label{eq:Gmap1}
\bx = G(\bxh) = \sum_{\bi\in \mathcal I} \bg_{\bi} \Bh_{\bi}(\bxh), 
\end{equation}
where $\bg_{\bi}$ are the spline (NURBS) control points, $\Bh_{\bi}$ are the basis functions, and $\bi$ is a $d$-dimensional multi-index, i.e., $\bi=(i_1,\dots, i_d)$. 

\begin{rem}
For simplicity, we consider only a single patch geometry map $G$ in this work. Multipatch parametrization techniques exist, we refer the reader for example to \cite{Kapl-2015-Isogeometric,Xu-2017-MultiPatch}.
\end{rem}

The basis functions possess the tensor product structure and we write
\begin{equation}\label{eq:Gmap2}
\Bh_{\bi}(\bxh) = \Bh_{i_1}(\xh_1) \dots \Bh_{i_d} (\xh_d),
\end{equation}
that are assumed, for the simplicity of the argument, to be piece-wise polynomial functions of the same degree $p$ in every variable. Let us denote by $\noe{EL}{i}$ 
the number of elements in the $i$-th parameter direction and define
\begin{equation}\label{eq:Gmap3}
\begin{array}{ccccc}
\Xi_i= ( 0 = & \underbrace{\xi_0, \dots,\xi_0,} & \xi_1, \xi_2, \dots, \xi_{\noe{EL}{i}-1}, 
& \underbrace{\xi_{\noe{EL}{i}},\dots,\xi_{\noe{EL}{i}}}& =1) \\
&  p+1 & & p+1 &
\end{array}
\end{equation}
the knot vector in the $i$-th direction. We assume that all the univariate splines have open knot vectors and are of the maximum continuity $C^{p-1}$, i.e., all the internal knots $\xi_1$, $\dots$, $\xi_{\noe{EL}{i}-1}$ are single knots. 

\begin{rem}
Conceptually, one can consider various degrees and continuities for univariate splines in each parameter direction. Assuming the same degree and continuity is not a limitation of the proposed method, it only simplifies the argument and offers a more convenient implementation.
\end{rem}
We denote by $\Sh_i$ the spline space of degree $p$ over a knot vector $\Xi_i$.
For its dimension we get
\begin{equation}\label{eq:DimUni}
\dim(\Sh_i)  = p+\noe{EL}{i} \quad i=1,\dots, d,
\end{equation}
and, due to the tensor product structure, the total number of degrees of freedom 
\begin{equation}\label{eq:DimMulti}
\noe{DOF}{} =  \prod_{i=1}^d \dim(S_i) 
\end{equation}
is the dimension of the approximate solution space
\begin{equation}\label{eq:PhysDiscr}
V_h = \textnormal{span}\lbrace \phi_{\bi}, \bi \in \mathcal I\rbrace \quad \textnormal{with} \quad \phi_{\bi}(\bx) =  \Bh_{\bi}( G^{-1}(\bx)), \quad \bx \in \Om.
\end{equation}
%
The Galerkin projection transforms the variational formulation \eqref{eq:Variational} into finding
\begin{equation}\label{eq:Galerkin1}
u_h \in V_h \quad \textnormal{such that} \quad a(u_h,v_h) = l(v_h) \quad \forall v \in V_h,
\end{equation}
where $h$ is the maximum element size. The elements of $V_h$ can be written as linear combinations of basis functions
\begin{equation}\label{eq:Galerkin2}
u_h = \sum_{\bi \in \mathcal{I}} u_{\bi}\phi_{\bi},
\end{equation}
with the coefficients vector $\bu=(u_1,\dots, u_{\noe{DOF}{}})$. 
We further introduce the pull-backs of the functions that are defined in the physical domain as
\begin{equation}\label{eq:PullBacks}
\uh_h (\bxh)= u_h(G(\bxh)) \quad \hat{f}(\bxh)= f(G(\bxh))
\quad \hat{\phi}(\bxh)= \phi(G(\bxh)), \quad \bxh \in \Omh
\end{equation}
and denote by $J$ the Jacobian matrix $J=\nh G$, $\nh$ being the gradient in the parameter domain.  Using the two forms \eqref{eq:Forms}, then the coefficients of  $u_h$ are a solution of the linear system
\begin{equation}\label{eq:LinSys}
\bK \bu = \bb,
\end{equation}
where $\bK$ is the \emph{stiffness matrix} 
\begin{equation}\label{eq:Stiffness}
\bK_{\bi\bj} = a(\phi_{\bi}, \phi_{\bj})  = \int_{\Omh} \nh \Bh_{\bi}(\bxh)\Transp J(\bxh)^{-1}  J(\bxh)^{\text{\sf\sc -t}}
\nh \Bh_{\bj}(\bxh) \, \vert\det(J(\bxh))\vert \de \bxh,
\end{equation}
and $\bb$ is the load vector
\begin{equation}\label{eq:Load}
\bb_{\bi} = l(\phi_{\bi})  = \int_{\Omh} \hat{f}  \Bh_{\bi}(\bxh) \,\vert\det(J(\bxh))\vert \de \bxh.
\end{equation}


For partial differential equations that contain zero order terms, the variational form contains also scalar products of basis functions that form the \emph{mass matrix}
\begin{equation}\label{eq:Mass}
\bM_{\bi\bj} = b(\phi_{\bi}, \phi_{\bj})  = \int_{\Omh} \Bh_{\bi}(\bxh) \Bh_{\bj}(\bxh) \,\vert\det(J(\bxh))\vert \de \bxh.
\end{equation}

\begin{figure}[!tb]
\vrule width0pt\hfill
    \begin{overpic}[width=.89\textwidth,angle=0]{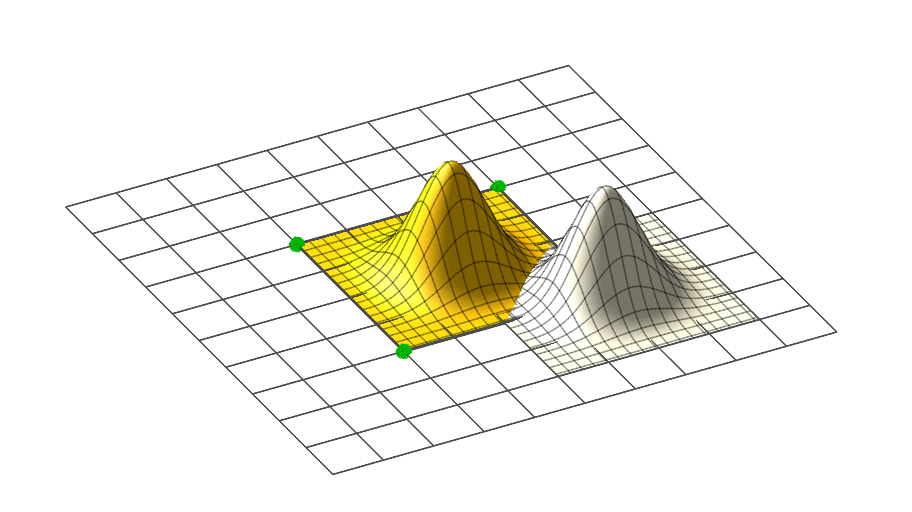}
        \put(42,25){\small$\Bh_{\bj}$}
        \put(61,22){\small$\Bh_{\bi}$}
        \put(8,38){\small$\RR^2$}
        \put(36,20){\small$\mathcal{H}$}
    \end{overpic}
\hfill \vrule width0pt\\[-2ex]
\Acaption{1ex}{Spline basis functions in the neighborhood of the weight function; $p=3$, $d=2$.
The 2D-grid defines a knot neighborhood (here uniform) of $(2p+1)^d$ basis functions $\Bh_{\bi}$ that have non-vanishing integrals with respect to the measure  $\Bh_{\bj}$, see \eqref{eq:Mass2}.}\label{fig:BasisWeight2D}
\end{figure}

\section{Integration via weighted quadrature}\label{sec:WQ}

We start our considerations with the numerical integration of the type 
\begin{equation}\label{eq:Mass2}
\int_{\Omt}  \Bh_{\bi}(\bxh) \Bh_{\bj}(\bxh) \, \de \bxh,
\end{equation}
which, due to the tensor product structure of $\Bh_{\bi}(\bxh)$, can be decomposed into a sequence of $d$ univariate integrations
\begin{equation}\label{eq:Decompo}
\int_0^1 \Bh_{i_1}(\xh_1) \Bh_{j_1}(\xh_1)  \Big[
\int_0^1 \Bh_{i_2}(\xh_2) \Bh_{j_2}(\xh_2)  \cdots  \Big[
\int_0^1 \Bh_{i_d}(\xh_d) \Bh_{j_d}(\xh_d) \de \xh_d \Big] \cdots \de \xh_2 \Big] \de \xh_1.
\end{equation}
We follow \cite{Sangali-2016-Weighted} where one basis function is the weight function for the integration. Due to the local properties of B-spline basis functions, the measure is positive only on the support of $\Bh_{\bj}$ and zero elsewhere. Defining $\mathcal{H} = \textnormal{supp}(\Bh_{\bj})$, $\mathcal{H} \subset \Omh$, a weighted Gaussian quadrature with respect to the measure $\Bh_{\bj}$ 
can be seen as a local quadrature mask acting only on the macroelement $\mathcal{H}$, see Fig.~\ref{fig:BasisWeight2D}.

\begin{figure}[!tb]
\vrule width0pt\hfill
    \begin{overpic}[width=.89\textwidth,angle=0]{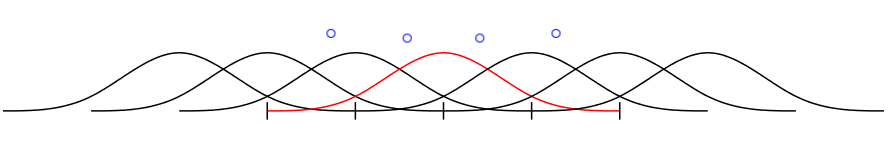}
        \put(48,12){\small$\Bh_{j}$}
        \put(8,10){\small$\Bh_{j-p}$}
        \put(85,10){\small$\Bh_{j+p}$}
        \put(30,0){\small $\xi_{j}$}
        \put(58,0){\small$\xi_{j+p}$}
        \put(68,0){\small$\xi_{j+p+1}$}
        \put(31,15){\small$[\tau_1,\omega_1]$}
    \end{overpic}
\hfill \vrule width0pt\\[-2ex]
\Acaption{1ex}{Weighted Gauss quadrature for univariate cardinal B-splines, $p=3$. Gaussian quadrature with respect to a measure $\Bh_j$ (red) requires $p+1$ quadrature points (blue dots) since only $2p+1$ basis functions $\Bh_{j-p}, \dots, \Bh_{j+p}$ have an overlapping support with the one of $\Bh_j$, $[\xi_j, \xi_j+p+1]$. 
 }\label{fig:BasisWeight}
 \end{figure}

The decomposition \eqref{eq:Decompo} allows us to consider a sequence of univariate integrals.
We aim to compute    
\begin{equation}\label{eq:WeightIntUni}
\int_{\textnormal{supp}(\Bh_{j})}  \Bh_{i}(\xh) \Bh_{j}(\xh) \, \de \xh
\end{equation}
by deriving a Gauss quadrature rule with respect to a non-negative measure $\mu = \Bh_{j}(\xh) \, \de \xh$. The Gauss must be exact for all B-spline functions $\Bh_{i}$ that have non-zero support on $\textnormal{supp}(\Bh_{j})$, see Fig.~\ref{fig:BasisWeight}.

The number of spline basis functions that have non-zero support on  $\textnormal{supp}(\Bh_{j})$
is $2p+1$ and therefore the weighted Gauss quadrature requires only $p+1$ quadrature points. Since $\Bh_{j}$ spans $p+1$ elements, this results in a one-node-per-element rule, regardless the degree. Therefore, our approach optimizes the scheme introduced in \cite{Sangali-2016-Weighted}, reducing the number of quadrature points in one parameter direction from $2p+1$ to $p+1$. Considering the tensor product structure of the weighted quadrature integration \eqref{eq:Decompo}, the total reduction ratio of the quadrature points is 
\begin{equation}\label{eq:Reduction}
\bigg(\frac{p+1}{2p+1}\bigg)^d.
\end{equation}

\begin{rem}
The reduction in \eqref{eq:Reduction} compares solely the number of quadrature points between our approach and \cite{Sangali-2016-Weighted}. Since that the quadrature points of \cite{Sangali-2016-Weighted} are the knots and midpoints, the evaluations of various basis functions at the quadrature points can be reused while our rule quadrature points differ for every weight function. Nonetheless, Section~\ref{sec:NUMER} shows that the number of evaluations is at least halved when using our approach. 
\end{rem}

We now derive specific rules for uniform $C^1$ quadratic and $C^2$ cubic isogeometric discretizations. With a slight abuse of notation, we omit the hat symbol over the basis functions and variables, remembering that all belong to the parameter domain.

\subsection{Weighted Gaussian quadrature -- Mass matrix terms}\label{ssec:WeightedGauss}

We seek an \emph{$m$-point }Gauss quadrature rule
\begin{equation}\label{eq:GaussQuadW}
\int_a^b f(x) w(x) \, \mathrm{d}x = \sum_{i=1}^{m} \omega_i f(\tau_i) + R_{m}(f),
\end{equation}
with respect to the non-negative weight function $w(x)$
such that the rule is exact for any function $f(x)$ from a spline space under consideration. 
The existence and uniqueness of such a rule has been investigated in \citep{Micchelli-1977}.
In our context, the weight function is represented by the spline basis function $\Bh_{j}$, the domain of integration is its support $[\xi_j, \xi_{j+p+1}]$, for some $j \in \ZZ$, and the number of Gaussian points is $m=p+1$. To meet the exactness constraint, we obtain a system 
\begin{equation}\label{eq:GaussQuadSys}
\int_{\xi_j}^{\xi_j+p+1} B_i(x) B_j(x) \, \mathrm{d}x = \sum_{k=1}^{p+1} \omega_k^j B_i(\tau_k^j)B_j(\tau_k^j) 
\quad \forall i = j-p, \dots, j+p,
\end{equation}
where the quadrature points $\tau_k^j$ and weights $\omega_k^j$ relate to the weight function $B_j$ and are the unknowns of an under-constrained system that consists of 
$2p+1$ equations and $2(p+1)$ unknowns.

\begin{rem}
According to  \cite[Theorem 3.1]{Micchelli-1977}, there exists a unique quadrature rule for splines with respect to a positive measure if the corresponding system is well-constrained (same number of unknowns and constraints). Our rules are Gaussian as they guarantee exactness with a minimum number of quadrature points, however, they are not unique as the rules are roots of underconstrained piece-wise polynomial systems.  
\end{rem}  

\subsubsection{$C^1$ quadratic elements}\label{ssec:MassMatrix}

\begin{figure}[!tb]
\vrule width0pt\hfill
    \begin{overpic}[width=.89\textwidth,angle=0]{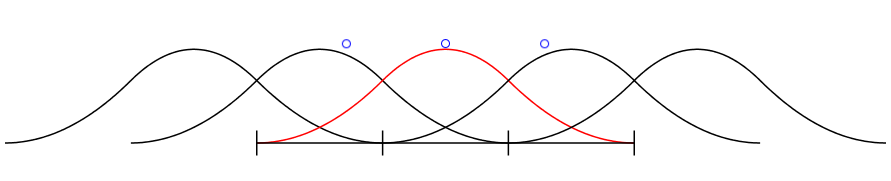}
        \put(46,11){\small$B_{3}$}
        \put(8,10){\small$B_{1}$}
        \put(88,10){\small$B_{5}$}
        \put(26,1){\small $0$}
        \put(72,1){\small$3$}
    \end{overpic}
\hfill \vrule width0pt\\[-2ex]
\Acaption{1ex}{Weighted Gaussian quadrature for $C^1$ quadratic cardinal spline ($p=2$); mass matrix term. The weight function $B_3$ (red) vanishes outside $[0,3]$ and a has non-trivial overlap with itself and four other basis functions. The Gaussian quadrature requires $p+1$ quadrature points (blue dots, \eqref{eq:WQquadratic}) that, due to symmetry, arise as a solution of the $(3\times 3)$ system \eqref{eq:QuadraticQuadrature}.}\label{fig:BasisWeight_p=2}
 \end{figure}

We start with the quadratic case ($p=2$). Consider the cardinal B-spline basis function $B_j$ with a non-zero support on $[0,3]$. There are $2p+1$ basis functions that interact with $B_j$ and therefore one needs three quadrature points to satisfy the conditions \eqref{eq:GaussQuadSys}, see Fig.~\ref{fig:BasisWeight_p=2}. There is one degree of freedom to choose either one quadrature point or weight. We impose the symmetry constraint and set $\tau_2=1.5$. Consequently, we solve a reduced, well-constrained $(3\times 3)$  system  
\begin{equation}\label{eq:QuadraticQuadrature}
\renewcommand{\arraystretch}{1.6}
\begin{array}{rcl}
\frac{1}{4}\omega_1 \tau_1^2 (1-\tau_1)^2 & = & \frac{1}{120},  \\
\frac{1}{4}\omega_1 \tau_1^2 (-2\tau_1^2  + 2\tau_1 + 1) +  \frac{3}{32} \omega_2  & = & \frac{13}{60},  \\
 \frac{1}{2}\omega_1\tau_1^4+\frac{9}{16}\omega_2 & =&\frac{11}{20} 
\end{array}
\end{equation}
with the unknowns $\tau_1$, $\omega_1$, and $\omega_2$.
Using computer algebra Maple,  \eqref{eq:QuadraticQuadrature} can be solved symbolically. Evaluating the solutions with a precision of twenty decimal digits, we obtain the weighted Gaussian rule
\begin{equation}\label{eq:WQquadratic}
\renewcommand{\arraystretch}{1.1}
\begin{array}{cclcccc}
\tau_1    & = & 0 .71241440095955149482, & & \omega_1 & = & 0.79410713110801847176, \\
\tau_2    & = & 1.5, &  & \omega_2 & = & 0.79595121334251753503.
\end{array}
\end{equation}

\subsubsection{$C^2$ cubic elments}\label{ssec:MassCubic}

For the cubic case, we proceed analogously.
The Gaussian quadrature for cubic cardinal B-splines on $[0,4]$ with a weight $B_j$ requires $p+1=4$ nodes, see Fig.~\ref{fig:BasisWeight}.  Due to symmetry, the system \eqref{eq:GaussQuadSys} can be reduced to a well-constrained $(4\times 4)$ system
\begin{equation}\label{eq:CubicQuadrature}
\renewcommand{\arraystretch}{1.6}
\begin{array}{rcl}
\frac{1}{36}\omega_1 \tau_1^3 (1-\tau_1)^3 & = & \frac{1}{5040},  \\
\omega_1 \tau_1^3 (\frac{1}{12} \tau_1^3 - \frac{1}{6} \tau_1^2 + \frac{1}{9})
-\frac{1}{12}\omega_2 (2-\tau_2)^3 (\tau_2^3 -4\tau_2^2 +4\tau_2 - \frac{4}{3})
  & = & \frac{1}{42},  \\
 -\frac{1}{12}\omega_1 \tau_1^3 (\tau_1^3 - \tau_1^2-\tau_1-\frac{1}{3})& -& \\
 \frac{1}{3}\omega_2( \tau_2^3 -\frac{9}{2}\tau_2^2 +6\tau_2-\frac{3}{2}) 
 (\tau_2^3 -4\tau_2^2 +4\tau_2 - \frac{4}{3}) & = & \frac{397}{1680},  \\
 \frac{1}{18}\omega_1\tau_1^6+\frac{1}{2}\omega_2
 (\tau_2^3 -4\tau_2^2 +4\tau_2 - \frac{4}{3})^2 & =&\frac{151}{315} 
\end{array}
\end{equation}
with the unknowns $\tau_1$, $\tau_2$, $\omega_1$, and $\omega_2$.
Solving it numerically gives 
\begin{equation}\label{eq:WQcubic}
\renewcommand{\arraystretch}{1.1}
\begin{array}{ccccccc}
\tau_1    & = & 0 .72289886179270511319, & & \omega_1 & = & 0.88863704203309628490, \\
\tau_2    & = & 1.58789880583487289415, &  & \omega_2 & = & 0.83494225417405959060.
\end{array}
\end{equation}

%
Let us point out that there are several difficulties with \eqref{eq:CubicQuadrature}. 
First of all, the system is built under the assumption that $\tau_1 \in [0,1]$ and $\tau_2 \in [1,2]$ which in general is not known. 
In this simple configuration with only four elements, one could follow the argumentation of \cite{Quadrature51-2017} and prove that there has to be a quadrature node in each element, but for higher degrees, or non-uniform knots, such an approach is not straightforward.  
But even more importantly, solving the system \eqref{eq:CubicQuadrature} numerically requires a good initial guess and, for example, with the initial guess 
\begin{equation}\label{eq:WQcubicIni}
\tau_1^{\textnormal{INI}}  = \frac{1}{3}, \quad  \tau_2^{\textnormal{INI}}  = \frac{5}{3},
\quad \omega_1^{\textnormal{INI}}  = \omega_2^{\textnormal{INI}}  =  1
\end{equation}
the Newton-Raphson finds a solution 
\begin{equation}\label{eq:WQcubicçWrong}
\renewcommand{\arraystretch}{1.1}
\begin{array}{ccccccc}
\tau_1    & = & 0.75698683155927590528, & & \omega_1 & = & 1.14740718959367949323, \\
\tau_2    & = & 2.30382606794266282352, & & \omega_2 & = & 0.74428414202245775486,
\end{array}
\end{equation}
that is not correct as it violates the assumption of $\tau_2 \in [1,2]$. For higher degrees, it is therefore obvious that a more elaborated strategy for finding the initial guess needs to be incorporated. This issue goes beyond the scope of the current paper and will be addressed in future work.

\subsection{Weighted Gaussian quadrature -- Stiffness matrix terms}\label{ssec:WeightedGauss}

To compute the stiffness matrix, one needs to integrate products of derivatives, see \eqref{eq:Stiffness}. Observe that derivatives of the spline basis functions change signs in their support. Therefore, to the best of our knowledge, there is no theoretical result that guarantees existence of the quadrature rule with respect to such a weight function \cite{Micchelli-1977}. However, we can build a weighted quadrature by using a piece-wise algebraic system.

\begin{figure}[!tb]
\vrule width0pt\hfill
    \begin{overpic}[width=.89\textwidth,angle=0]{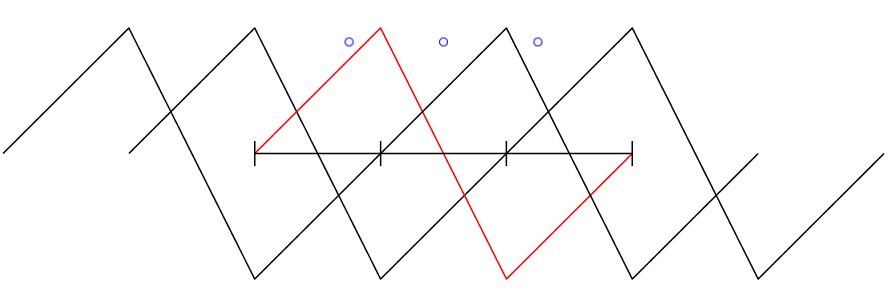}
        \put(44,29){\small$B'_{3}$}
        \put(2,16){\small$B'_{1}$}
        \put(92,14){\small$B'_{5}$}
        \put(28,12){\small $0$}
        \put(70,12){\small$3$}
    \end{overpic}
\hfill \vrule width0pt\\[-4ex]
\Acaption{1ex}{Weighted Gauss quadrature for $C^1$ quadratic cardinal spline ($p=2$); stiffness matrix terms. The piece-wise linear derivatives of the quadratic spline basis functions that interact with the weight function $B'_3$ (red) are shown. The weighted Gaussian quadrature (blue dots, \eqref{eq:WQquadraticStiff}) arises as a solution of the system \eqref{eq:QuadraticQuadratureStiff}.}\label{fig:BasisWeightStiffness_p=2}
 \end{figure}

The derivatives of the basis functions of degree $p$ still span $p+1$ elements,
and there are, same to the mass matrix case, $2p+1$ basis functions that interact with $B'_j$.
However, this spline space is only $2p$-dimensional, that is, the derivatives that have non-zero support on the support of $B'_j$ are linearly dependent. Having only $p$ quadrature points is not possible due to the symmetry and the fact that the first $p+1$ derivatives are linearly independent on the support of $B'_j$, see Fig.~\ref{fig:BasisWeightStiffness_p=2}. Therefore to obtain our sought weighted Gaussian rule, we build the exactness constraints as 
\begin{equation}\label{eq:GaussQuadSysStiff}
\int_{\xi_j}^{\xi_j+p+1} B'_i(x) B'_j(x) \, \mathrm{d}x = \sum_{k=1}^{p+1} \omega_k^j B'_i(\tau_k^j) B'_j(\tau_k^j) 
\quad \forall i = j-p, \dots, j+p,
\end{equation}
that contains $2p+2$ unknowns and $2p+1$ constraints, only $2p$ being independent. Since the support of $B'_j$ is $p+1$ elements, we postulate, that for uniform knots, there is one node in every element and build the system \eqref{eq:GaussQuadSysStiff} accordingly. 

\subsubsection{$C^1$ quadratic elements}\label{ssec:StiffnessMatrix}

If the initial space is $C^1$ quadratic, the derivatives are piece-wise linear and due to symmetry, the second quadrature point is the middle point $\tau_2=\frac{3}{2}$, see Fig.~\ref{fig:BasisWeightStiffness_p=2}. The system \eqref{eq:GaussQuadSysStiff} becomes simply
\begin{equation}\label{eq:QuadraticQuadratureStiff}
\renewcommand{\arraystretch}{1.6}
\begin{array}{rcl}
\omega_1 \tau_1 (1-\tau_1)  & = & \frac{1}{6},  \\
\omega_1 \tau_1 (1-2\tau_1)  & = & \frac{1}{3},  \\
 (\omega_1 + \omega_2) \tau_1^2 & =& 1, 
\end{array}
\end{equation}
with the unknowns $\tau_1$, $\omega_1$, and $\omega_2$. The system 
gets solved in radicals with a unique solution
\begin{equation}\label{eq:WQquadraticStiff}
\renewcommand{\arraystretch}{1.6}
\begin{array}{cclcccl}
\tau_1    & = & \frac{3}{4}, & & \omega_1 & = & \frac{8}{9}, \\
\tau_2    & = & \frac{3}{2}, & & \omega_2 & = & \frac{8}{9}.
\end{array}
\end{equation}

\subsubsection{$C^2$ cubic elements}\label{ssec:StiffnessMatrix}

In the cubic case, we have four quadrature points with  seven constraints 
in \eqref{eq:GaussQuadSysStiff}. Due to symmetry, we build a $(4\times 4)$ system 
\begin{equation}\label{eq:CubicQuadratureStiff}
\renewcommand{\arraystretch}{1.6}
\begin{array}{rcl}
\frac{1}{4}\omega_1 \tau_1^2 (1-\tau_1)^2 & = & \frac{1}{120},  \\
\frac{1}{4}\omega_1 \tau_1^3 (4-3\tau_1) - \frac{3}{4}\omega_2(\tau_2-2)^3(\tau_2-\frac{2}{3}) & = & \frac{1}{5},  \\
\frac{1}{4}\omega_1 \tau_1^2 (3\tau_2^2-2\tau_2-1) +3\omega_2(\tau_2-1)(\tau_2-2)^2(\tau_2-\frac{2}{3}) & = & \frac{1}{8},  \\
 \frac{1}{2}\omega_1\tau_1^4+\frac{9}{2}\omega_2(\tau_2-2)^2 (\tau_2-\frac{2}{3})^2& =&\frac{2}{3}, 
\end{array}
\end{equation}
with the unknowns $\tau_1$, $\tau_2$, $\omega_1$, and $\omega_2$. 
The system admits one parameter family of solutions, $\tau_1$ being a root of
\begin{equation}\label{eq:CubicStiffTau1W}
30x^4\omega_1-60x^3\omega_1+30x^2\omega_1-1 = 0,
\end{equation}
while $\omega_1$ being the free parameter. Setting $\omega_1=1$, we obtain 
two roots of \eqref{eq:CubicStiffTau1W} that lie outside $[0,1]$, while the other two 
\begin{equation}\label{eq:CubicStiffTau1}
\tau_1 = \frac{1}{2} \pm \frac{1}{30} \sqrt{225-30\sqrt{30}}
\end{equation}
are admissible solutions inside $[0,1]$. These solutions
define two possible quadrature rules for this specific choice of $\omega_1$. Finally, setting $\tau_1$ as the smaller root of \eqref{eq:CubicStiffTau1}, the quadrature rule for stiffness matrix reads as
\begin{equation}\label{eq:WQcubicStiff}
\renewcommand{\arraystretch}{1.1}
\begin{array}{cclcccl}
\tau_1    & = & 0.24033518882038592858 & & \omega_1 & = & 1 \\
\tau_2    & = & 1.16015740029939774803 & & \omega_2 & = & 0.86030876544418464920.
\end{array}
\end{equation}

\begin{figure}[!tb]
\vrule width0pt\hfill
    \begin{overpic}[width=.89\textwidth,angle=0]{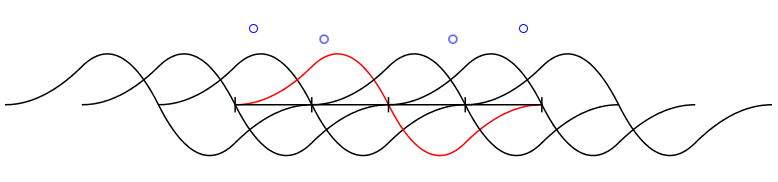}
        \put(44,17){\small$B'_{4}$}
        \put(2,12){\small$B'_{1}$}
        \put(96,6){\small$B'_{7}$}
        \put(28,7){\small $0$}
        \put(71,10){\small$4$}
    \end{overpic}
\hfill \vrule width0pt\\[-4ex]
\Acaption{1ex}{Weighted Gauss quadrature for $C^2$ cubic splines ($p=3$); stiffness matrix terms. Blue dots represent the quadrature rule \eqref{eq:WQcubicStiff} that comes from the system \eqref{eq:CubicQuadratureStiff}.}\label{fig:BasisWeightStiffness_p=3}
 \end{figure}

\begin{rem}
We derive weighted Gaussian quadrature rules only for internal basis functions, i.e., when the weight functions are not affected by the boundary. The number of degrees of freedom with support on the boundary is negligible for large meshes. For boundary elements, the corresponding system would need to be modified accordingly, or one can use standard Gauss rules there.
\end{rem}

\begin{rem}
Our rules are exact for affine geometric mappings as for them the Jacobians in \eqref{eq:Stiffness} and \eqref{eq:Mass} are constants. For general, non-affine maps, our rules compute only approximate values that correspond to a spline approximation of rational functions.
\end{rem}

%

\section{Numerical examples}\label{sec:NUMER}

In this section, we present numerical examples that demonstrate the efficiency and validate the accuracy of the proposed quadrature rules. Our rules exactly integrate the spline spaces associated to the mass and stiffness matrices and therefore integrate exactly (up to machine precision) all the matrix entries.

\begin{figure}[!tb]
   	\begin{overpic}[width=1.0\textwidth,angle=0]{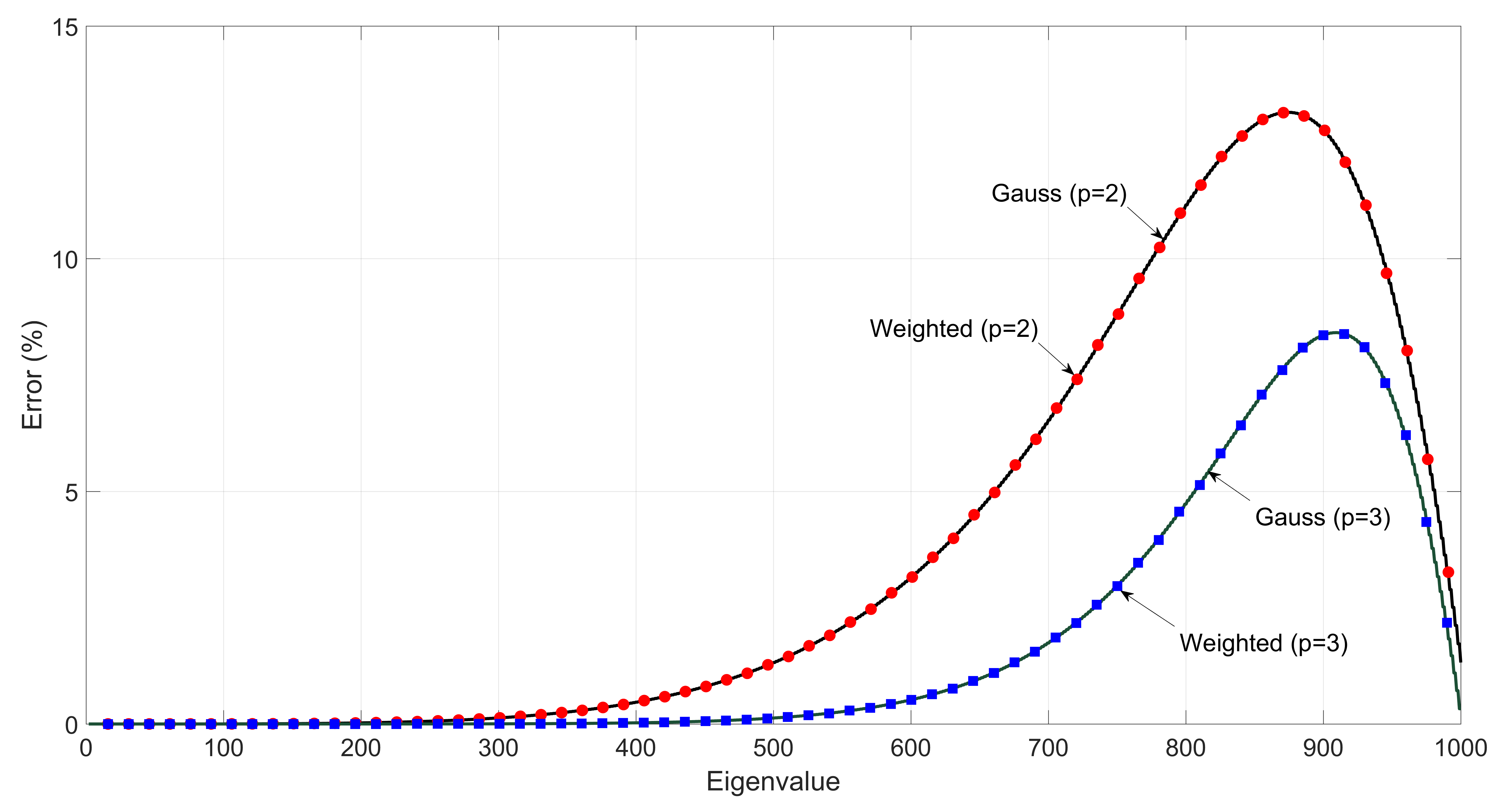}
    \end{overpic}
\Acaption{1ex}{Approximation errors of the quadratic and cubic isogeometric discretizations using the standard Gaussian quadratures (lines) and our weighted quadrature rules (markers). The quadrature rules described in \citep{Sangali-2016-Weighted} leads to the same results and are not shown here for brevity.}\label{fig:ResultComparison}
\end{figure}

In the following plots, we show the eigenvalue approximation errors that are fundamental for error estimation in many boundary- and initial-value problems \citep{Cottrell-2009}.
Figure~\ref{fig:ResultComparison} compares the approximation errors of the standard $C^1$ quadratic and $C^2$ cubic isogeometric elements when the mass and stiffness matrices are assembled by the standard Gauss quadrature rule and the weighted Gaussian quadratures proposed in the previous section. In this example, we consider a one-dimensional elliptic eigenvalue problem discretized on a uniform mesh with 1000 elements. The example numerically validates the exactness of our rules as it shows almost identical error when compared to the standard Gauss integration. The maximum absolute difference in the terms of the mass and stiffness matrices created by both quadratures is of order $10^{-15}$.

\begin{figure}[!tb]
	\begin{subfigure}[t]{0.5\textwidth}
    	\begin{overpic}[width=1.0\textwidth,angle=0]{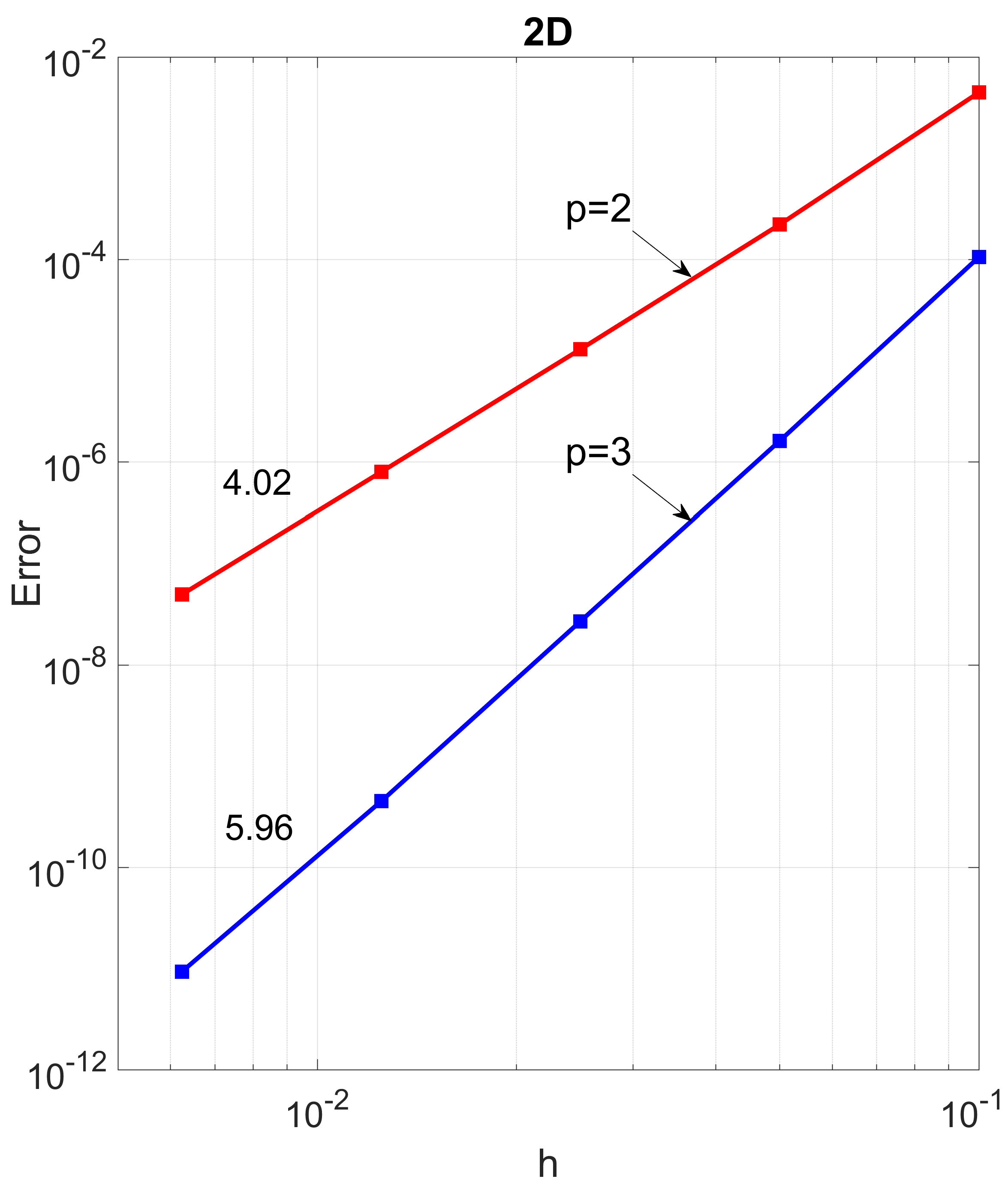}
	    \end{overpic}
	\end{subfigure}
    \hfill
    \begin{subfigure}[t]{0.5\textwidth}
       	\begin{overpic}[width=1.0\textwidth,angle=0]{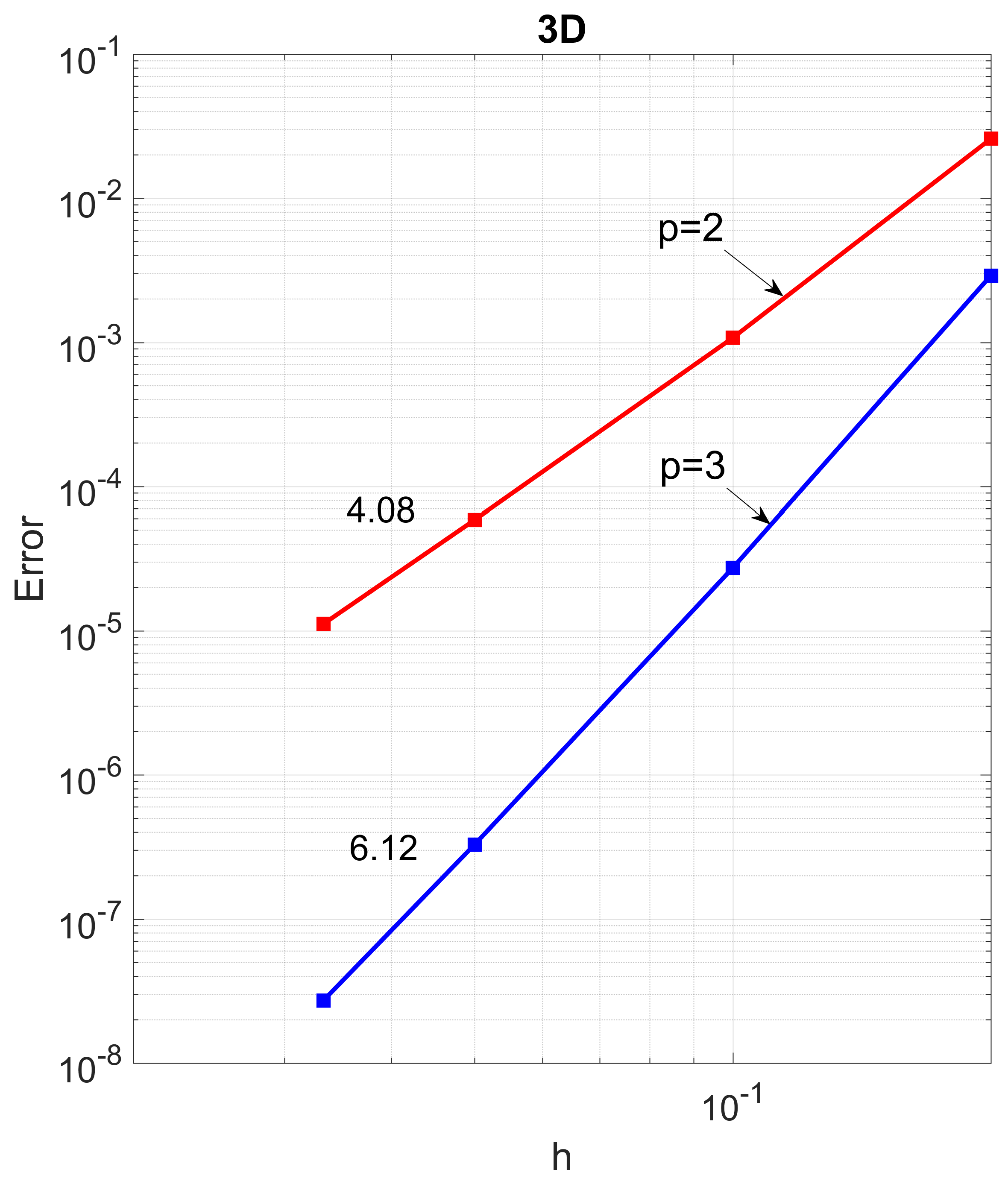}
	    \end{overpic}
	\end{subfigure}
\Acaption{1ex}{Convergence of the $10$-th eigenvalue of the 2D (left) and 3D (right) problems using the quadratic (red line) and cubic (blue line) isogeometric elements. The mass and stiffness matrices were assembled using the weighted quadrature rules introduced in Section~\ref{sec:WQ}.}\label{fig:Result2DConvergence}
\end{figure}

Figure~\ref{fig:Result2DConvergence} shows the convergence of the 2D and 3D eigenvalue problems. We apply the weighted quadrature rules for mass and stiffness matrices for $C^1$ quadratic (\eqref{eq:WQquadratic} and \eqref{eq:WQquadraticStiff}) and $C^2$ cubic (\eqref{eq:WQcubic} and \eqref{eq:WQcubicStiff})
isogeometric discretizations. The convergence rates obtained numerically are close to the theoretical order of $2p$ for the fully-integrated case \citep{deng2017quadratures}.

\begin{figure}[!tb]
	\centering
   	\begin{overpic}[width=0.5\textwidth,angle=0]{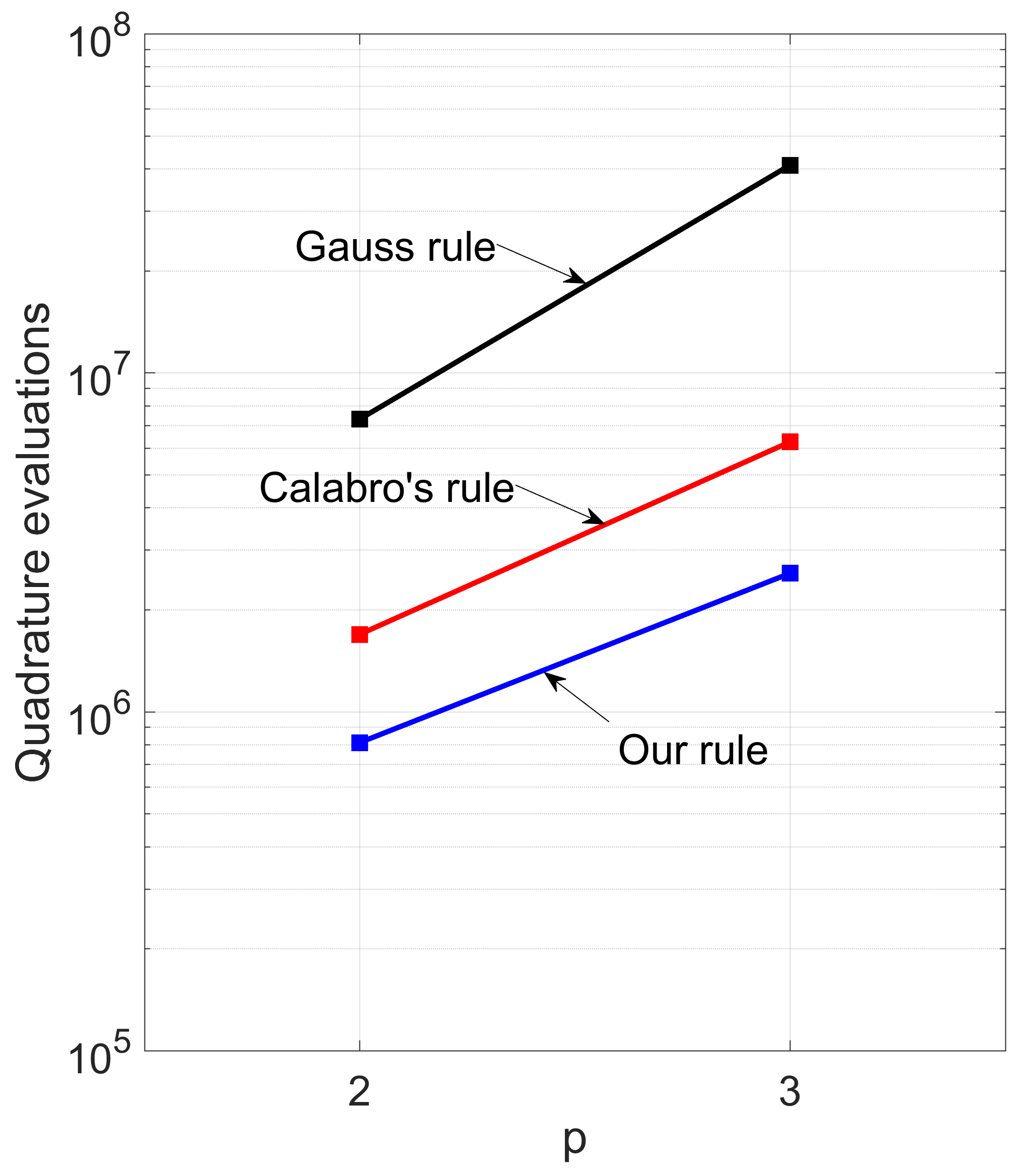}
    \end{overpic}
\Acaption{1ex}{Number of quadrature evaluations for a test 2D problem discretized on a 100x100 mesh: standard Gaussian quadrature (black line), the quadrature rule of \citep{Sangali-2016-Weighted} (red line) and our weighted quadrature rule (blue line).}\label{fig:ResultQuadCalls}
\end{figure}

In order to show the computational gain when using the weighted Gaussian quadratures, we compare the number of basis function evaluations at the quadrature points required to build the mass matrix. Figure~\ref{fig:ResultQuadCalls} shows the total number of evaluations for a test 2D problem discretized on a 100x100 mesh. Our weighted quadratures outperform both the standard element-wise Gaussian algorithm and also \cite{Sangali-2016-Weighted}. The computational gain for the stiffness matrix is similar. Our rules require less than half the quadrature evaluations compared to \citep{Sangali-2016-Weighted} for both quadratic and cubic elements. This numerical result is in accordance with the theoretical estimate \eqref{eq:Reduction}. 

\begin{rem}
The exact comparison values are $2.08$ and $2.44$, which are slightly worse than the theoretical estimates $\frac{25}{9}$ and $\frac{49}{16}$ obtained from \eqref{eq:Reduction}. The reason is that the rules of Calabr\'{o} et al. \citep{Sangali-2016-Weighted} use the knots and the middle of the elements as quadrature points and therefore, in this uniform setup, some function evaluations may vanish. To make a fair comparison, we excluded from our counting the cases where the function evaluation is redundant. 
For example, a scenario where the weight function vanishes at the middle knot, see Fig.~\ref{fig:BasisWeightStiffness_p=3}.
\end{rem}


\section{Conclusion}\label{sec:con} %

We present weighted Gaussian quadrature rules for  the mass and stiffness matrix assembly of $C^1$ quadratic and $C^2$ cubic spline discretizations over uniform knot vectors. Our rules are the solutions of a set of nonlinear piece-wise polynomial systems. For quadratic elements, the rules factorize to a closed form solution, while for cubics we find solutions numerically. Our weighted quadrature rules
require only one quadrature point per element and therefore further reduce the cost of mass and stiffness matrix assembly when compared to \citep{Sangali-2016-Weighted} which requires two quadrature points per element. 

As a future work, we aim to focus on higher degree elements. Already in the cubic case, the numerical solver did not find a correct solution from a reasonable initial guess and therefore a more elaborated strategy to initialize the numerical solver seems to be unavoidable. Another direction goes along the analysis conducted in \cite{Puzyrev-2017-Dispersion, deng2017quadratures, deng2017quadratures2} and the fact that the exact integration of the mass matrix is not needed to achieve the optimal convergence rate. Therefore, weighted rules that underintegrate mass terms but use even fewer quadrature points are worth of further investigation. 

%
%

\section*{Acknowledgements}

The first author has been partially supported by the Basque Government through the BERC 2014-2017 program, by Spanish Ministry of Economy and Competitiveness MINECO: BCAM Severo Ochoa excellence accreditation SEV-2013-0323, and the Project of the Spanish Ministry of Economy and Competitiveness with reference MTM2016-76329-R (AEI/FEDER, EU). 
This publication was made possible in part by the CSIRO Professorial Chair in Computational Geoscience at Curtin University and the Deep Earth Imaging Enterprise Future Science Platforms of the Commonwealth Scientific Industrial Research Organisation, CSIRO, of Australia. Additional support was provided by the European Union's Horizon 2020 Research and Innovation Program of the Marie Skłodowska-Curie grant agreement No. 644202.


\bibliographystyle{elsarticle-harv}\biboptions{square,sort,comma,numbers}
\bibliography{Quadrature}

\end{document}